\documentclass[twoside,leqno,10pt]{article}
\usepackage{graphics,ifthen}
\usepackage{latexsym}
\usepackage{amsmath}
\usepackage{amssymb}
\usepackage{mathrsfs}

\begin{document}
\input{latexP.sty}
\input{referencesP.sty}
\input epsf.sty

\def\ind{\stackrel{\mathrm{ind}}{\sim}}
\def\iid{\stackrel{\mathrm{iid}}{\sim}}

\def\Definition{\stepcounter{definitionN}\
    \Demo{Definition\hskip\smallindent\thedefinitionN}}
\def\EndDefinition{\EndDemo}
\def\Example#1{\Demo{Example {\rm #1}}}
\def\EndExample{\qed\EndDemo}
\def\Category#1{\centerline{\Heading #1}\rm}
\
\def\e{\text{\hskip1.5pt e}}
\newcommand{\eps}{\epsilon}
\newcommand{\proofs}{\noindent {\bf Proof of Theorem 2:\ }}
\newcommand{\prove}{\noindent {\bf Proofs of points (ii):\ }}
\newcommand{\remarks}{\noindent {\bf Remarks:\ }}
\newcommand{\note}{\noindent {\bf Note:\ }}
\newcommand{\examp}{\noindent {\bf Example:\ }}
\newcommand{\Lower}[2]{\smash{\lower #1 \hbox{#2}}}
\newcommand{\ben}{\begin{enumerate}}
\newcommand{\een}{\end{enumerate}}
\newcommand{\bi}{\begin{itemize}}
\newcommand{\ei}{\end{itemize}}
\newcommand{\hp}{\hspace{.2in}}

\numberwithin{equation}{section}
\newtheorem{thm}{Theorem}
\newtheorem{defin}{Definition}[section]
\newtheorem{prop}{Proposition}[section]
\newtheorem{lem}{Lemma}[section]
\newtheorem{cor}{Corollary}
\newcommand{\rb}[1]{\raisebox{1.5ex}[0pt]{#1}}
\newcommand{\mc}{\multicolumn}
\newcommand{\Acr}{\mathscr{A}}
\newcommand{\Bcr}{\mathscr{B}}
\newcommand{\Ucr}{\mathscr{U}}
\newcommand{\Gcr}{\mathscr{G}}
\newcommand{\Dcr}{\mathscr{D}}
\newcommand{\CS}{\mathscr{C}}
\newcommand{\Fcr}{\mathscr{F}}
\newcommand{\Icr}{\mathscr{I}}
\newcommand{\Lcr}{\mathscr{L}}
\newcommand{\Mcr}{\mathscr{M}}
\newcommand{\ncr}{\mathscr{n}}
\newcommand{\Ncr}{\mathscr{N}}
\newcommand{\Ocr}{\mathscr{O}}
\newcommand{\Pcr}{\mathscr{P}}
\newcommand{\Qcr}{\mathscr{Q}}
\newcommand{\Scr}{\mathscr{S}}
\newcommand{\Tcr}{\mathscr{T}}
\newcommand{\Xcr}{\mathscr{X}}
\newcommand{\Vcr}{\mathscr{V}}
\newcommand{\Ycr}{\mathscr{Y}}
\newcommand{\qcr}{\mathscr{q}}
\newcommand{\scr}{\mathscr{s}}
\newcommand{\indic}{\mathbb{I}}
\newcommand{\E}{\mathbb{E}}
\newcommand{\F}{\mathbb{F}}
\newcommand{\I}{\mathbb{I}}
\newcommand{\Q}{\mathbb{Q}}
\newcommand{\X}{\mathbb{X}}
\newcommand{\Pe}{\mathbb{P}}
\newcommand{\M}{\mathbb{M}}
\newcommand{\Wbb}{\mathbb{W}}

\def\Beta{\text{Beta}}
\def\Dir{\text{Dirichlet}}
\def\DP{\text{DP}}
\def\P{{\bf p}}
\def\fhat{\widehat{f}}
\def\GA{\text{gamma}}
\def\ind{\stackrel{\mathrm{ind}}{\sim}}
\def\iid{\stackrel{\mathrm{iid}}{\sim}}
\def\J{{\bf J}}
\def\K{{\bf K}}
\def\min{\text{min}}
\def\N{\text{N}}
\def\p{{\bf p}}
\def\U{{\bf U}}
\def\W{{\bf W}}
\def\S{{\bf S}}
\def\T{{\bf T}}
\def\y{{\bf y}}
\def\t{{\bf t}}
\def\m{{\bf m}}
\def\X{{\bf X}}
\def\Y{{\bf Y}}
\def\tps{\mbox{\scriptsize ${\theta H}$}}   
\def\ups{\mbox{\scriptsize ${P_{\theta}(g)}$}}   
\def\vps{\mbox{\scriptsize ${\theta}$}}   
\def\vups{\mbox{\scriptsize ${\theta >0}$}}   
\def\hps{\mbox{\scriptsize ${H}$}}   
\def\rps{\mbox{\scriptsize ${(\theta+1/2,\theta+1/2)}$}}   
\def\sps{\mbox{\scriptsize ${(1/2,1/2)}$}}   

\newcommand{\reals}{{\rm I\!R}}
\newcommand{\PR}{{\rm I\!P}}
\def\Z{{\bf Z}}
\def\yy{{\mathcal Y}}
\def\rr{{\mathcal R}}
\def\BP{\text{beta}}
\def\ts{\tilde{t}}
\def\js{\tilde{J}}
\def\gs{\tilde{g}}
\def\fs{\tilde{f}}
\def\ys{\tilde{Y}}
\def\ps{\tilde{\mathcal {P}}}

\def\Report{Lancelot F. James and Marc Yor}
\def\Author{}
\pagestyle{myheadings}
\markboth{\Author}{\Report}
\thispagestyle{empty}

\bct\Heading  Tilted stable subordinators, Gamma time changes and
Occupation Time of rays by Bessel Spiders \lbk\lbk\smc Lancelot F.
James and Marc Yor \footnote{\eightit \rm Supported in
part by grants HIA05/06.BM03 and DAG04/05.BM56 of the HKSAR.
\eightit Corresponding authors address.
                \rm The Hong Kong University of Science and Technology,
Department of Information and Systems Management, Clear Water Bay,
Kowloon, Hong Kong.
\rm lancelot\at ust.hk\\
\indent\eightit Keywords and phrases.
                \rm
Bessel spiders, Gamma process, Tilted stable subordinators.
          }
\lbk\lbk \BigSlant The Hong Kong University of Science and
Technology\\and Universit\'e Pierre et Marie Curie- Paris VI.\rm \lbk 
\ect \Quote We exhibit, in the form of some identities in law,
some connections between tilted stable subordinators, time-changed
by independent Gamma processes and the occupation times of Bessel
spiders, or their bridges. These identities in law are then
explained thanks to excursion theory. \EndQuote
\rm
\section{Introduction} \subsection{Aim of this work} The genesis of this paper is our interest in
a class of subordinators $(T_{t},t\ge 0)$ which enjoy the
following properties: (i) $(T_{t},t\ge 0)$ is of GGC (Generalized
Gamma convolution) type. That is, it has no drift and its L\'evy
measure is of the form \begin{equation} \label{levyden}\theta
\frac{dx}{x}\E[{\mbox e}^{-xG}] \end{equation} for $G$, a random
variable such that $\E[{(\log(1/G))}^{+}]<\infty,$ and $\theta>0;$
(ii) each marginal law of $T_{t}$, for fixed $t,$ can be described
explicitly. Of course, this second demand is not mathematically
very precise but we are seeking descriptions such as: $$ T_{t}\sim
X_{t}Y_{t}, {\mbox { for fixed t }} $$ where $X_{t}$ and $Y_{t}$
and have classical (e.g. beta, gamma) distributions. [Throughout,
we shall write: $U\sim V$, to mean that $U$ and $V$ are
identically distributed.] As we see above the law of the
subordinator $(T_{t})$, which can be written more precisely as
$(T^{(\theta,G)}_{t})$ depends both on $G$ and $\theta.$ However,
we note that $(T^{(\theta,G)}_{t}, t\ge 0)$ is equivalent (in law)
to $(T^{(1,G)}_{\theta t}, t\ge 0).$

\Remark The general class of GGC random variables has been
extensively studied by Thorin~(1977) and Bondesson~(1992). These
are infinitely divisible random variables and hence one can
naturally construct subordinators from them. We note that in
general the L\'evy measure of such GGC subordinators may be
defined as $x^{-1}\int_{0}^{\infty}{\mbox e}^{-xy}U(dy)$, where
$U$ is a sigma-finite measure. That is to say it does not
necessarily correspond to, or arise from, a random variable. The
case where $U(\infty)=\infty$, contains many important classes of
random variables, including positive stable random variables of
index $0<\alpha<1$ and the GIG (Generalized Inverse Gaussian)
class. See Bondesson~(1992) for many properties and examples of
GGC random variables. \EndRemark
\subsection{First examples} Here, we next give two examples of
subordinators which satisfy the desired properties (i) and (ii)
mentioned above. First let us take $\theta =1/2$, then:
\begin{example} $ G:=G_{1/2}\sim \beta_{(1/2,1/2)},$ where
$\beta_{(a,b)}$ denotes a beta random variable, with parameters
$(a,b).$\end{example}\begin{example} $G:=1/G_{1/2}\sim
1/\beta_{(1/2,1/2)}.$\end{example} We denote these subordinators
by $\varepsilon_{1}$ and $\varepsilon_{2}$, respectively~(where
here $\varepsilon$ stands for example). Then,  as we shall see in
section 3, there are the equations; \Eq
\label{eq1.2}\E[\exp(-\lambda
\varepsilon_{1}(t))]={(\sqrt{1+\lambda}-\sqrt{\lambda})}^{t}=\frac{1}{{(\sqrt{1+\lambda}+\sqrt{\lambda})}^{t}}
\EndEq and
$\varepsilon_{1}(t)\sim\gamma_{t/2}/\beta_{(1/2,(1+t)/2)}.$ \Eq
\label{eq1.3} \E[\exp(-\lambda
\varepsilon_{2}(t))]=\frac{2^{t}}{{(1+\sqrt{1+\lambda})}^{t}}
\EndEq and
$\varepsilon_{1}(t)\sim\gamma_{t/2}\beta_{((1+t)/2,(1+t)/2)}.$
\Remark Interestingly, the processes $\varepsilon_{1}$ and
$\varepsilon_{2},$ have appeared previously. $\varepsilon_{1}$ is
a special case of the random variables recently studied in
Bertoin, Fujita, Roynette, and Yor~(2006) and Fujita and
Yor~(2006). The fact that
$\varepsilon_{1}(t)\sim\gamma_{t/2}/\beta_{(1/2,(1+t)/2)},$ as
indicated in~\mref{eq1.2}, is demonstrated in Roynette and
Yor~(2006) where connections to Feller~(1966) are also noted. The
fact that
$\varepsilon_{2}(t)\sim\gamma_{t/2}\beta_{((1+t)/2,(1+t)/2)},$ as
indicated in~\mref{eq1.3}, follows from a result of Cifarelli and
Melilli~(2000), concerning Dirichlet mean functionals, as
discussed in James~(2006).\EndRemark
\subsection{A general class of variables G}
The fact that $\beta_{(1/2,1/2)}$ is the arcsine law, and that
this law is that of the time spent in $\reals_{+}$ (or
$\reals_{-}$) by real valued Brownian motion, led us to look for
generalizations of the above examples, where the variable $G$
in~\mref{levyden} would be replaced by a variable distributed as
the time spent in $\reals_{+}$ by a symmetrized Bessel process and
even to look for further generalizations. As shown in Barlow,
Pitman and Yor~(1989), the time spent in $\reals_{+}$ by a
symmetrized Bessel process with dimension $d=2(1-\alpha)$, for
$0<\alpha<1$, is distributed as;
$$
\frac{{\mathbf T}_{\alpha}}{{\mathbf T}_{\alpha}+{\mathbf
T}'_{\alpha}},
$$
where $\textbf{T}_{\alpha}$ and ${\mathbf T}'_{\alpha}$ are two
independent standard stable $(\alpha)$ random variables, i.e.
$\E[\exp(-\lambda {\mathbf T}_{\alpha})]=\exp(-\lambda^{\alpha}).$
More generally,  we looked for a description of the subordinators
satisfying~\mref{levyden} with \begin{equation} \label{genG}
G\overset{d}=\frac{\sum_{j=1}^{N}\nu_{j}{\mathbf
T}^{(j)}_{\alpha}}{\sum_{j=1}^{N}\mu_{j}{\mathbf
T}^{(j)}_{\alpha}} \end{equation} where all the coefficients are
assumed to be non-negative, and the $({\mathbf T}^{(j)}_{\alpha},
j=1,2,\ldots, N)$ are $N$ iid standard stable $(\alpha)$
variables.

It turns out that the subordinator $T_{t}$ associated with $G$
given in~\mref{genG} may be obtained from $$
\CS_{u}:=\CS^{(\mu,\nu)}_{u}=\sum_{i=1}^{N}\mu_{i}T^{(\nu_{i})}_{\alpha}(u)
$$ where the $(T^{(\nu_{i})}_{\alpha}(u);u\ge 0)$ are constructed
from $N$ independent stable $(\alpha)$ subordinators, each of them
being Esscher transformed with parameter $(\nu_{i})$, that is,
$$
\E[\exp(-\lambda
T^{(\nu_{i})}_{\alpha}(u))]=\exp(-u[{(\nu_{i}+\lambda)}^{\alpha}-\nu^{\alpha}_{i}]).$$
In fact, $(T_{t})$ is constructed as $(\CS_{\gamma_{t}}, t\ge 0)$
where $(\gamma_{t}, t\ge 0)$ is a standard Gamma process
independent of $(\CS_{u}, u\ge 0).$ We call the process
$(\CS_{\gamma_{t}}, t\ge 0)$ a (generalized) positive Linnik
process with index $\alpha,$ and parameters $(\mu,\nu).$ The name
we use comes from its link to the class of random variables
expressible as
$\textbf{L}_{\alpha,t}=\gamma^{1/\alpha}_{t}{\mathbf T}_{\alpha},$
defined for each fixed $t>0,$ where $\gamma_{t}$ and ${\mathbf
T}_{\alpha}$ are independent. Naturally these random variables are
equivalent in distribution to a positive stable $(\alpha)$
subordinator time changed by a gamma process $\gamma_{t},$ for
each fixed $t.$ The random variables, $\textbf{L}_{\alpha,t}$, are
known in the literature[see for instance Bondesson~(1992, p.38)
and Devroye~(1990, 1996)] and are sometimes called generalized
positive Linnik random variables. These random variables have a
host of interesting properties and distributional representations.
Huillet~(2000) discusses some results related to the subordinator.
James~(2006) provides an extensive recent study of
$\textbf{L}_{\alpha, t}$ and its Esscher transformed version which
connects with the work of Barlow, Pitman and Yor~(1989), and our
present exposition. In this sense, $\textbf{L}_{\alpha, t}$ and
its Esscher transformed version arise as special cases of
$(\CS_{\gamma_{t}}, t\ge 0).$ Moreover, $(\CS_{\gamma_{t}}, t\ge
0),$ also contains an interesting subclass of the models discussed
in the recent related work of Bertoin, Fujita, Roynette and
Yor~(2006).  We shall show in Section 4 that we are able to
describe precisely the one-dimensional marginals of the
generalized positive Linnik process, $(\CS_{\gamma_{t}}, t\ge 0),$
which therefore answers positively our demands (i) and (ii) made
in (1.1). These one-dimensional marginals are closely related with
the occupation times of the Bessel spiders living in webs
constituted by $N$ rays as described in Barlow, Pitman and
Yor~(1989). Once this was noticed, it became of interest to
explain this coincidence. This is what we do in the last two
sections of the paper.
\subsection{Organization of the paper} We now describe the organization of the remainder
of this paper: in Section 2, we recall notation and results about
Bessel spiders and their occupation times of rays; in Section 3,
we show that we can describe the one-dimensional marginals and the
L\'evy measures of the generalized Linnik process, in Section 4,
we compare the results obtained in Section 2 and 3; in the final
Section 5, we give an explanation,  based on excursion theory, for
the close relationship between the distributions we have
encountered in Sections 2 and 3.
\section{On the occupation times of rays by a Bessel spider}
In this section, we simply recall the definition of a Bessel
spider, with index $\alpha\in (0,1)$, and whose state space is
$E,$ the union of $N$ half lines, or rays, $I_{1},\ldots, I_{N},$
originating from a common point which we denote as $0.$ In fact,
these random processes are indexed by both $\alpha\in (0,1), $ and
$\textbf{p}=(p_{1},\ldots, p_{N}),$ a generic probability on
$\{1,2,\ldots, N\}.$ Now, we may introduce $\Scr(\alpha,
p)=(S_{t},t\ge 0)$, a Bessel spider with index $\alpha,$ and
probability $p$, which lives on $E$, as, informally, a process
$(S_{t},t\ge 0)$ which behaves like a BES$(\alpha)$ process on
each of the rays, $I_{1},\ldots, I_{N},$ and which, when coming
back to $0$, rides off in all directions at once, and chooses its
ray $I_{i}$, with probability $p_{i}.$ Of course, this is only a
heuristic description as the point $0$ is regular for itself (with
respect to the Markov process $(S_{t})$), but rigorous
constructions have been provided by Barlow, Pitman and Yor~(1989),
with the help in particular of excursion theory.

The case $\alpha=1/2$ is particularly interesting, as the
$(1/2,\textbf{p})$ spider behaves like a Brownian motion in each
of its rays.; We also note that for $\alpha=1/2$, and $N=2,$
$\p=(p_{1},p_{2})$, $\Scr(1/2,\p)$ may be identified with the skew
Brownian motion of, say, parameter $p_{1}$ (or $p_{2}$!).

We also recall that $|S_{t}|:=\textrm{dist}(S_{t},0)$ is a
\textsc{BES}$(\alpha)$ process (in particular, for $\alpha=1/2, $
it is a reflecting Brownian motion), and $(S_{t}, t\ge 0)$ admits
a local time at $0,$ $(L_{t}, t\ge0)$ such that
$({|S_{t}|}^{2\alpha}-L_{t}, t\ge 0)$ is a martingale; we refer
the reader to Donati-Martin, Roynette, Vallois, and Yor~(2005) for
a discussion of different choices of this local time found in the
literature. In the following theorem we describe the joint laws of
: $\{A^{(i)}_{t}=\int_{0}^{t}ds \indic(S_{s}\in I_{i}); 1\leq
i\leq N; L_{t}\}$, for fixed $t, $ as well as the sequence
$\{a^{(i)}_{t}; i\leq N, \ell_{t}\},$ the same quantity, but now
considered for the standard spider's bridge.
\begin{thm}\begin{enumerate}\item[(1)]There is the identity in law
$$
\{A^{i}_{1}, 1\leq i\leq
N\}\overset{law}=\{\frac{p_{i}^{1/\alpha}\textbf{T}^{(i)}_{\alpha}}{\sum_{j}p^{1/\alpha}_{j}\textbf{T}^{(j)}_{\alpha}},
i\leq N;
\frac{1}{\sum_{j}p^{1/\alpha}_{j}\textbf{T}^{(j)}_{\alpha}}\},
$$
where $\E[(\exp(-\lambda
\textbf{T}^{(j)}_{\alpha})]=\exp(-\lambda^{\alpha})$
\item[(2)]The law of $\{a^{(i)}_{1}; i\leq N,
\ell^{1/\alpha}_{1}\}$ is absolutely continuous w.r.t that of;
$\{(A^{i}_{1}, i\leq N); L^{1/\alpha}_{1}\}$ with density
$\Gamma(1-\alpha)L_{1}$. This may be written as: $\forall$
$f\ge0,$ Borel ,$:[0,1]^{N}\times\mathbb{R}_{+}\rightarrow
\mathbb{R}_{+},$
$$
\E[f((a^{(i)},i\leq N); \ell^{1/\alpha}_{1})]=
\E[f(\frac{p^{1/\alpha}_{i}\textbf{T}^{(i)}_{\alpha}}{\sum_{j}p_{j}^{1/\alpha}\textbf{T}^{(j)}_{\alpha}};
\i\leq N;
\frac{1}{\sum_{j}p^{1/\alpha}_{j}\textbf{T}^{(j)}_{\alpha}})\frac{\Gamma(1+\alpha)}{\sum_{j}p^{1/\alpha}\textbf{T}^{(j)}_{\alpha}}].
$$

\end{enumerate}
\end{thm}
In the next section, we shall see, in a different context, similar
quantities-in particular, the ratios: $
{\textbf{T}^{(i)}_{\alpha}}/{\sum_{j}\nu_{j}\textbf{T}^{(j)}_{\alpha}}$
appear.
\section{An identity about multivariate positive generalized Linnik processes} In this
section, we consider
$$
\CS^{(\nu)}_{\alpha}(u)=(T^{(\nu_{i})}_{\alpha}(u); 1\leq i\leq
N), u\ge 0,
$$
an $\mathbb{R}^{N}_{+}-$valued process, whose components are
independent $\nu_{i}$-Esscher transforms of
stable$(\alpha)$-subordinators. Let, furthermore, $(\gamma_{t},
t\ge 0)$ denote a standard gamma process independent of
$\CS^{(\nu)}_{\alpha}.$ The next theorem offers a description of
the joint law of  $(\CS^{(\nu)}_{\alpha}(\gamma_{t}), \gamma_{t}),
$ in terms of a sequence of iid stable-$(\alpha)$ variables,
$(\textbf{T}^{(i)}_{\alpha}, i\leq N)$, and an independent gamma
$(t\alpha)$ variable.
\begin{thm}\begin{enumerate}\item[1)]For any
$F:\mathbb{R}^{N}_{+}\rightarrow \mathbb{R}_{+},$ and
$g:\mathbb{R}_{+}\rightarrow \mathbb{R}_{+},$ Borel, one has:
\begin{equation} \label{eq3.1}
\E[F((\CS^{(\nu)}_{\alpha}(\gamma_{t}))g(\gamma_{t})]=\frac{\Gamma(1+\alpha
t)}{\Gamma(1+t)}\E\left[F\left((\frac{\gamma_{t}}{\nu
\cdot\textbf{T}_{\alpha}})\textbf{T}_{\alpha}\right)g({(\frac{\gamma_{t}}{\nu
\cdot\textbf{T}_{\alpha}})}^{\alpha})\frac{\exp({(\frac{\gamma_{t}}{\nu
\cdot\textbf{T}_{\alpha}})}^{\alpha}(\sigma-1))}{{(\nu\cdot
\textbf{T}_{\alpha})}^{t\alpha}}\right]
\end{equation}
where $\sigma=\sum_{i=1}^{N}\nu^{\alpha}_{i},$ and
$\nu\cdot\textbf{T}_{\alpha}=\sum_{j=1}^{N}\nu_{j}\textbf{T}^{(j)}_{\alpha}.$
In particular for $\sigma=1$, the formula simplifies as:
\begin{equation} \label{eq3.2}
\E[F((\CS^{(\nu)}_{\alpha}(\gamma_{t}))g(\gamma_{t})]=\frac{\Gamma(1+\alpha
t)}{\Gamma(1+t)}\E\left[F\left((\frac{\gamma_{t}}{\nu
\cdot\textbf{T}_{\alpha}})\textbf{T}_{\alpha}\right)g({(\frac{\gamma_{t}}{\nu
\cdot\textbf{T}_{\alpha}})}^{\alpha})\frac{1}{{(\nu\cdot
\textbf{T}_{\alpha})}^{t\alpha}}\right].
\end{equation}
\item[2)]If $\Lcr(dx_{1},\ldots, dx_{N})$ denotes the L\'evy measure
of the $N$-dimensional L\'evy process
$(\CS^{(\nu)}_{\alpha}(\gamma_{t}),t\ge 0)$ then for any
$F:\mathbb{R}^{N}_{+}\rightarrow \mathbb{R}_{+}$ continuous, with
compact support in ${(0,\infty)}^{N}$:
$$
\langle\Lcr,F\rangle=\lim_{t\downarrow
0}(\frac{1}{t}\E[F(\CS^{(\nu)}_{\alpha}(\gamma_{t}))])=\alpha\int_{0}^{\infty}\frac{du}{u}{\mbox
e}^{-u}\E[F(u\frac{\textbf{T}_{\alpha}}{\nu\cdot
\textbf{T}_{\alpha}})]
$$
\end{enumerate} \end{thm}
\Remark Note that the last formula follows from
$$
\E[F(X_{t})]=\E[\int_{0}^{t}ds LF(X_{s})],
$$
where $(X_{t})$ is a L\'evy process, starting at $0$, $F$ a
regular function, equal to $0$ on a neighborhood of $0,$ and $L$
is the infinitesimal generator of $X$. If $X$ has no drift, and no
Brownian component, then
$$
LF(0)=\langle\Lcr,F\rangle.
$$
\EndRemark

\begin{cor}
In the case $\sigma=1,$ the variables
$\nu\cdot\CS^{(\nu)}_{\alpha}(\gamma_{t})$ and
$\CS^{(\nu)}_{\alpha}(\gamma_{t})/{(\gamma_{t})}^{1/\alpha}$ are
independent and satisfy:
\begin{enumerate}
\item[(i)]$\nu\cdot
\CS^{(\nu)}_{\alpha}(\gamma_{t})\overset{d}=\gamma_{\alpha t}.$
\item[(ii)]\begin{equation}
\label{eq3.3}
\E\left[H\left(\frac{1}{\gamma^{1/\alpha}_{t}}\CS^{(\nu)}_{\alpha}(\gamma_{t})\right)\right]=\frac{\Gamma(1+\alpha
t)}{\Gamma(1+t)}\E\left[H(\textbf{T}_{\alpha})\frac{1}{{(\nu\cdot
\textbf{T}_{\alpha})}^{t\alpha}}\right],\end{equation} for any
positive measurable function $H:\mathbb{R}^{N}_{+}\rightarrow
\mathbb{R}_{+},$ or equivalently:
\begin{equation}
\label{eq3.4}
\frac{1}{\gamma^{1/\alpha}_{t}}\CS^{(\nu)}_{\alpha}(\gamma_{t})\overset{d}=\CS^{(\nu
\gamma^{1/\alpha}_{t})}_{\alpha}(1),
\end{equation}
where on the RHS of \mref{eq3.4}, the law of the variable
appearing there is the mixture, with respect to the law of $(\nu
\gamma^{1/\alpha}_{t})$ of the laws of $\CS^{(\mu)}_{\alpha}(1)$,
i.e. the $\mu$-Esscher transform of $\CS^{(0)}_{\alpha}(1).$
$\gamma_{t}$ is assumed to be independent of the process:
$(\CS^{(\mu)}_{\alpha}(1), \mu\ge 0).$
\end{enumerate}
\end{cor}
\begin{cor} Let $\mu_{i}\ge 0,$ $\nu_{i}>0;$ $1\leq i\leq N.$ The generalized positive Linnik process
$$\CS^{(\mu,\nu)}_{\alpha}(\gamma_{t})\equiv\sum_{i=1}^{N}\mu_{i}T^{(\nu_{i})}_{\alpha}(\gamma_{t})\equiv\mu\cdot
\CS^{(\nu)}_{\alpha}(\gamma_{t}))$$ is a GGC subordinator with
L\'evy measure:
$$
\alpha x^{-1}\E\left[\exp\left(-x \frac{(\nu\cdot
\textbf{T})}{(\mu\cdot \textbf{T})}\right)\right]dx.$$
\end{cor}
 Corollary 2 follows immediately from the expression of $\Lcr$ (or
 rather $\langle\Lcr, F\rangle$) in point 2) of Theorem 2. Likewise, Corollary
 1 is easily deduced from formula~\mref{eq3.2} in theorem 2, by
 easy changes of variables. We now give a proof of theorem
 2:\\
 \begin{proofs} Of course,~\mref{eq3.2} follows immediately
 from~\mref{eq3.1}, so it suffices to prove~\mref{eq3.1}. We first
 use the absolute continuity relationship between the laws of
 $(\CS^{(\nu)}_{\alpha}(s); s\leq u)$ and $(\CS^{(0)}_{\alpha}(s); s\leq
 u),$ which yields that the LHS of ~\mref{eq3.1} equals:
 $$
\E[F((\CS^{(0)}_{\alpha}(\gamma_{t}),\gamma_{t})\exp\{-\sum_{i=1}^{N}(\nu_{i}T^{(0)}_{\alpha}(\gamma_{t})-\gamma_{t}\nu^{\alpha}_{i}\}]=
\E[F(\gamma^{1/\alpha}_{t}\textbf{T}_{\alpha},\gamma_{t})\exp(-\gamma^{1/\alpha}_{t}(\nu\cdot\textbf{T}_{\alpha})+\sigma
\gamma_{t})].
$$
(by scaling and independence),
$$
=\frac{1}{\Gamma(t)}\int_{0}^{\infty}dx x^{t-1}{\mbox
e}^{x(\sigma-1)}\E[F(x^{1/\alpha}\textbf{T}_{\alpha},
x)\exp(-x^{1/\alpha}(\nu\cdot \textbf{T}_{\alpha}))]
$$
(then, changing variables $x=y^{\alpha}$),
$$
=\frac{\alpha}{\Gamma(t)}\int_{0}^{\infty}dy
y^{\alpha-1}{(y^{\alpha})}^{t-1}{\mbox
e}^{y^{\alpha}(\sigma-1)}\E[F(y\textbf{T}_{\alpha},
y^{\alpha})\exp(-y(\nu\cdot \textbf{T}_{\alpha}))]$$ (then
changing variables: $y=z/(\nu\cdot \textbf{T}_{\alpha})$),
$$
=\frac{\alpha \Gamma(\alpha
t)}{\Gamma(t)}\E\left[F(\frac{\gamma_{\alpha t}}{\nu\cdot
\textbf{T}_{\alpha}}\textbf{T}_{\alpha},{(\frac{\gamma_{\alpha
t}}{\nu\cdot
\textbf{T}_{\alpha}})}^{\alpha})\frac{\exp({(\frac{\gamma_{\alpha
t}}{\nu
\cdot\textbf{T}_{\alpha}})}^{\alpha}(\sigma-1))}{{(\nu\cdot
\textbf{T}_{\alpha})}^{\alpha t}}\right]
$$
which yields the desired result.
 \end{proofs}

 \begin{remark}
 We also note that the preceding argument yields, but we do not
 give the details, the following $1$-parameter extension of
 formula~\mref{eq3.1}
\begin{equation} \label{eq3.5}
\E[F((\CS^{(\nu)}_{\alpha}(\frac{1}{m}\gamma_{t}))g(\gamma_{t})]=\frac{\Gamma(1+\alpha
t)}{\Gamma(1+t)}\E\left[m^{t}F\left((\frac{\gamma_{\alpha t}}{\nu
\cdot\textbf{T}_{\alpha}})\textbf{T}_{\alpha}\right)g(m{(\frac{\gamma_{\alpha
t}}{\nu
\cdot\textbf{T}_{\alpha}})}^{\alpha})\frac{\exp({(\frac{\gamma_{\alpha
t}}{\nu
\cdot\textbf{T}_{\alpha}})}^{\alpha}(\sigma-m))}{{(\nu\cdot
\textbf{T}_{\alpha})}^{\alpha t}}\right]
\end{equation}
\end{remark}
\subsection{Examples}
We now illustrate Theorem 2 and Corollary 2 in the particular case
$N=2$ and $\alpha=1/2$, for which the law of $\textbf{T}_{1/2}$ is
explicit (recall $\E[{\mbox e}^{-\lambda \textbf{T}_{1/2}}]={\mbox
e}^{-\sqrt{\lambda}}$). We shall also assume
$\sqrt{\nu_{1}}+\sqrt{\nu_{2}}=1.$ \Example{1} Set $\nu_{1}=1,$
$\nu_{2}=0,$ $\mu_{1}=\mu_{2}=1;$ $\E[\exp(-\lambda
T_{t})]={(\sqrt{1+\lambda}-\sqrt{\lambda})}^{t}.$
\begin{prop}\begin{enumerate}\item[(i)]The L\'evy measure of
 the process $(T_{t})$ is
 $$
 \frac{1}{2}\frac{du}{u}\E[\exp(-u\beta_{(1/2,1/2)})]
 $$
\item[(ii)]The distribution of the variable $T_{t}$ is that of
$$
\frac{\gamma_{t/2}}{\beta_{(\frac{1}{2},\frac{1+t}{2})}}.
$$
\end{enumerate}
\end{prop}
\EndExample \Example{2}Set $\nu_{1}=\nu_{2}=1/4,$ $\mu_{1}=1/4,$
$\mu_{2}=0;$ $\E[\exp(-\lambda
T_{t})]={(2/(1+\sqrt{1+\lambda}))}^{t}.$
\begin{prop}\begin{enumerate}\item[(i)]The L\'evy measure of
 the process $(T_{t})$ is
 $$
 \frac{1}{2}\frac{du}{u}\E[\exp(-u/\beta_{(1/2,1/2)})]
 $$
\item[(ii)]The distribution of the variable $T_{t}$ is that of
$$
{\gamma_{t/2}}{\beta_{(\frac{1+t}{2},\frac{1+t}{2})}}.$$
\end{enumerate}

\end{prop}
\EndExample
\begin{example} (This is the general case for $\alpha=1/2, N=2$): Set $\sqrt{\nu_{1}}+\sqrt{\nu_{2}}=1,$ $\mu_{1},\mu_{2}\ge 0;$

\begin{prop}\begin{enumerate}\item[(i)]The L\'evy measure of
 the process $(T_{t})$ is
 $$
 \frac{1}{2}\frac{du}{u}\E[\exp(-u(\frac{\nu_{1}\textbf{T}_{1/2}+\nu_{2}\textbf{T}'_{1/2}}{\mu_{1}\textbf{T}_{1/2}+\mu_{2}\textbf{T}'_{1/2}}))]
 $$
\item[(ii)]The law of the variable $T_{t}$ satisfies: $\forall f\ge 0$
\begin{equation}
\label{eq3.6}
\E[f(T_{t})]=C_{t}\E\left[f\left(\gamma_{t/2}\frac{\mu_{1}\hat{\gamma}'_{(1+t)/2}+\mu_{2}\hat{\gamma}_{(1+t)/2}}
{\nu_{1}\hat{\gamma}'_{(1+t)/2}+\nu_{2}\hat{\gamma}_{(1+t)/2}}\right)\frac{1}
{{(\nu_{1}\hat{\gamma}'_{\frac{(1+t)}{2}}+\nu_{2}\hat{\gamma}_{\frac{(1+t)}{2}})}^{t/2}}\right
] \end{equation} for some universal constant $C_{t},$ and$\gamma,$
$\hat{\gamma},$ $\hat{\gamma}'$ are independent gamma processes.

\end{enumerate}
\end{prop}
\end{example}
Rather than proving Propositions 3.1, 3.2, 3.3 in that order, we
shall first prove Proposition 3.3, and then see that the results
of  Proposition 3.1 and 3.2 follow. In fact, we only need to
prove the point (ii) in each of the Propositions.\\
\begin{prove}\begin{itemize}\item[(a)] From Theorem 2, we know:
$$
\E[f(T_{t})]=\frac{\Gamma(1+\frac{t}{2})}{\Gamma(1+t)}\E\left[f\left(\gamma_{t/2}\frac{\mu_{1}\textbf{T}^{(1)}_{1/2}+\mu_{2}\textbf{T}^{(2)}_{1/2}}
{\nu_{1}\textbf{T}^{(1)}_{1/2}+\nu_{2}\textbf{T}^{(2)}_{1/2}}\right)\frac{1}
{{(\nu_{1}\textbf{T}^{(1)}_{1/2}+\nu_{2}\textbf{T}^{(2)}_{1/2})}^{t/2}}\right
].
$$
Since $\textbf{T}^{(i)}_{1/2}\overset{d}=c/\gamma_{1/2}$ for
$i=1,2$, we get, with obvious notation:
$$
\E[f(T_{t})]=C_{t}\E\left[f\left(\gamma_{t/2}\frac{\mu_{1}\hat{\gamma}'_{1/2}+\mu_{2}\hat{\gamma}_{1/2}}
{\nu_{1}\hat{\gamma}'_{1/2}+\nu_{2}\hat{\gamma}_{1/2}}\right){\left(\frac{\hat{\gamma}_{\frac{1}{2}}\hat{\gamma}'_{\frac{1}{2}}}
{{(\nu_{1}\hat{\gamma}'_{\frac{1}{2}}+\nu_{2}\hat{\gamma}_{\frac{1}{2}})}}\right)}^{t/2}\right
]
$$
which simplifies to~\mref{eq3.6}.
\item[(b)]In the case of Proposition 3.2, we have:
$\nu_{1}=\nu_{2},$ and then
$\hat{\gamma}'_{{(1+t)}/{2}}/(\hat{\gamma}'_{{(1+t)}/{2}}+\hat{\gamma}_{{(1+t)}/{2}})$
is a \textsc{beta}((1+t)/2,(1+t)/2) random variable, independent
from $(\hat{\gamma}'_{{(1+t)}/{2}}+\hat{\gamma}_{{(1+t)}/{2}}).$
\item[(c)]In the case of Proposition 3.1, the formula~\mref{eq3.6} in point
(ii), Proposition 3.3., simplifies as:
$$
\E[f(T_{t})]=C_{t}\E\left[f\left(\gamma_{t/2}\frac{\hat{\gamma}'_{(1+t)/2}+\hat{\gamma}_{(1+t)/2}}
{\hat{\gamma}'_{(1+t)/2}}\right)\frac{1}
{{(\hat{\gamma}'_{\frac{(1+t)}{2}})}^{t/2}}\right ]
$$
$$
=\E\left[f\left(\gamma_{t/2}\frac{\hat{\gamma}'_{1/2}+\hat{\gamma}_{(1+t)/2}}
{\hat{\gamma}'_{1/2}}\right)\right ]
$$
$$
=\E\left[f\left(\gamma_{t/2}\frac{1}
{\beta_{(\frac{1}{2},\frac{1+t}{2})}}\right)\right ].
$$
\end{itemize}
\end{prove}
It might also be of interest to discuss the case $N=2$ and general
$0<\alpha<1$ in the same vein, because the law of
$(\textbf{T}_{\alpha}/\textbf{T}'_{\alpha})$ is simple, as shown
by Lamperti~(1958). Let us denote by $\textbf{R}_{\alpha}$ this
ratio. Then we have the following:
\begin{prop}($N=2$, $0<\alpha<1,$
$\nu^{\alpha}_{1}+\nu^{\alpha}_{2}=1$).
\begin{enumerate}\item[(i)]The L\'evy measure of the process is:
$$
\alpha\frac{du}{u}\E\left[\exp\left(-u\frac{\nu_{1}\textbf{R}_{\alpha}+\nu_{2}}{\mu_{1}\textbf{R}_{\alpha}+\mu_{2}}\right)\right]
$$
\item[(ii)]The law of $T_{t}$ satisfies $\forall f\ge0,$
\begin{equation}
\label{eq3.7}
\E[f(T_{t})]=C_{t}\E\left[f\left(\gamma_{t/2}\frac{\mu_{1}\textbf{R}_{\alpha}+\mu_{2}}
{\nu_{1}\textbf{R}_{\alpha}+\nu_{2}}\right)\frac{1}
{{(\nu_{1}\textbf{R}_{\alpha}+\nu_{2})}^{\alpha
t}}\frac{1}{({\textbf{T}'_{\alpha})}^{\alpha t}}\right ],
\end{equation}
\end{enumerate}
where, on the RHS, $\gamma_{t/2}$ is independent from the
pair~$(\textbf{T}_{\alpha},\textbf{T}'_{\alpha})$ and
$\textbf{R}_{\alpha}=\textbf{T}_{\alpha}/\textbf{T}'_{\alpha}.$
\end{prop}
\subsection{Looking for a simplification of Proposition 3.4}Using a
conditional expectation argument we can replace
$(\textbf{T}'_{\alpha})^{-\alpha t}$ appearing in~\mref{eq3.7}
with
$$
h_{\alpha
t}(\textbf{R}_{\alpha}):=\E[{(\textbf{T}'_{\alpha})}^{-\alpha
t}|\textbf{R}_{\alpha}].$$ Naturally this expression has utility
only if it has a tractable form. Here we obtain such an expression
by using Kanter's~(1975) explicit representation of the density of
a positive stable random variable, in conjunction with Lamperti's
expression for the density of $\textbf{R}_{\alpha}$. This form of
the stable density is apparently not well-known. First let
$f_{\alpha}$ denote the density of a unilateral stable-$(\alpha)$
density. Then generically the conditional density of
$\textbf{T}'_{\alpha}|\textbf{R}_{\alpha}=r$ is given by
$$
f_{\textbf{T}'_{\alpha}|\textbf{R}_{\alpha}}(s|r)=\frac{sf_{\alpha}(rs)f_{\alpha}(s)}{f_{\textbf{R}_{\alpha}}(r)}
$$
Now setting
$$
K_{\alpha}(u)={\left(\frac{\sin(\pi\alpha u)}{\sin(\pi
u)}\right)}^{-\frac{1}{1-\alpha}}{\left(\frac{\sin((1-\alpha) \pi
u)}{\sin(\pi \alpha u)}\right)},
$$
it follows from Kanter~(1975)[see also Devroye~(1996)] that
$$
f_{\alpha}(s)=\frac{\alpha}{1-\alpha}s^{-1/(1-\alpha)}\int_{0}^{1}{\mbox
e}^{-{s}^{-\frac{\alpha}{1-\alpha}}K_{\alpha}(u)}K_{\alpha}(u)du.
$$
That is
$\textbf{T}'_{\alpha}\overset{d}={(K_{\alpha}(U)/\textbf{e})}^{(1-\alpha)/\alpha},$
where $U$ is a Uniform$[0,1]$ random variable independent of
$\textbf{e},$ which is exponential$(1).$ Additionally from
Lamperti~(1958)[see, also Chaumont and Yor~(2003, p. 116)] we
obtain,
$$
f_{\textbf{R}_{\alpha}}(r)=\frac{\sin(\pi
\alpha)}{\pi}\frac{r^{\alpha-1}}{r^{2\alpha}+2r^{\alpha}\cos(\pi
\alpha)+1}{\mbox { for }}r>0.
$$
Combining these points we arrive at the following result
\begin{prop}\begin{enumerate} \item[(i)]The conditional density of
$\textbf{T}'_{\alpha}|\textbf{R}_{\alpha}=r$ is expressible as
$$
\frac{\pi r^{-\frac{\alpha(2-\alpha)}{1-\alpha}}}{\sin(\pi
\alpha)}{\left(\frac{\alpha}{1-\alpha}\right)}^{2}s^{-\frac{2}{1-\alpha}+1}\int_{0}^{1}\int_{0}^{1}{\mbox
e}^{-{s}^{-\frac{\alpha}{1-\alpha}}C_{\alpha}(u_{1},u_{2}|r)}K_{\alpha}(u_{1})K_{\alpha}(u_{2})du_{1}du_{2}
$$
where
$C_{\alpha}(u_{1},u_{2}|r)=[r^{-\frac{\alpha}{1-\alpha}}K_{\alpha}(u_{1})+K_{\alpha}(u_{2})]$
\item[(ii)]It follows that
$h_{\alpha
t}(\textbf{R}_{\alpha}):=\E[{(\textbf{T}'_{\alpha})}^{-\alpha
t}|\textbf{R}_{\alpha}]$ equals, \begin{equation} \label{eq3.8}
 h_{\alpha t}(\textbf{R}_{\alpha})=\frac{\pi
{(\textbf{R}_{\alpha})}^{-\frac{\alpha(2-\alpha)}{1-\alpha}}}{\sin(\pi
\alpha)}{\left(\frac{\alpha\Gamma(t(1-\alpha)+2)}{1-\alpha}\right)}
D_{\alpha,\alpha t}(\textbf{R}_{\alpha})
\end{equation}
where $D_{\alpha,\alpha t}(\textbf{R}_{\alpha})
=\int_{0}^{1}\int_{0}^{1}{[C_{\alpha}(u_{1},u_{2}|\textbf{R}_{\alpha})]}^{-(t(1-\alpha)+2)}K_{\alpha}(u_{1})K_{\alpha}(u_{2})du_{1}du_{2}.
$
\end{enumerate}
\end{prop}
This leads to an alternative characterization of $(T_{t})$ in
proposition 3.4 as follows;
\begin{prop}The distribution of the variable $T_{t}$ in Proposition 3.4 with $N=2$, $0<\alpha<1,$
$\nu^{\alpha}_{1}+\nu^{\alpha}_{2}=1,$ satisfies $\forall f\ge0,$
$$
\E[f(T_{t})]=C_{t}\E\left[f\left(\gamma_{t/2}\frac{\mu_{1}\textbf{R}_{\alpha}+\mu_{2}}
{\nu_{1}\textbf{R}_{\alpha}+\nu_{2}}\right)\frac{h_{\alpha
t}(\textbf{R}_{\alpha})}
{{(\nu_{1}\textbf{R}_{\alpha}+\nu_{2})}^{\alpha t}}\right ]
$$
where $h_{\alpha t}(\textbf{R}_{\alpha})$ is given
by~\mref{eq3.8}.
\end{prop}
\section{A comparison of Theorems 1 and 2}
This comparison is only partial, as here we shall use the
hypothesis: $\sigma=\sum_{i=1}^{N}\nu^{\alpha}_{i}=1,$ so that we
can associate with the Esscher sequence $(\nu_{i})_{i\leq N}$ the
probabilities $(p_{i})_{i\leq N}$ defined by
$p_{i}=\nu^{\alpha}_{i}.$ Then a simple cross-inspection of
Theorems 1 and 2 shows that, taking $t=1$,
\begin{equation}\label{eq4.1}
\{(\nu_{i}\CS^{(\nu_{i})}_{\alpha}(\textbf{e});i\leq
N),\textbf{e}\}\overset{d}=\{(\gamma_{\alpha}a^{(i)}_{1}; i\leq
N), {(\gamma_{\alpha})}^{\alpha}\ell_{1}\},
\end{equation}
where, on the LHS, \textbf{e} is an independent $\exp(1)$
variable, and on the RHS, $\gamma_{\alpha}$ is independent of the
Spider's Bridge. We shall now transform the identity in
law~\mref{eq4.1}, or rather its RHS, by using the two following
facts:
\begin{enumerate}
\item[(i)]if $g=\sup\{s<1:S_{s}=0\},$ then
$g\overset{d}=\beta_{(\alpha,1-\alpha)}$,
\item[(ii)]$\gamma_{\alpha}=\beta_{(\alpha,1-\alpha)}\textbf{e}$
where $\beta_{(\alpha,1-\alpha)}$ denotes a beta variable with
parameters~$(\alpha,1-\alpha),$ which is independent from
\textbf{e}.
\end{enumerate}
Thus the identity in law~\mref{eq4.1} may be rephrased as:
\begin{equation}\label{eq4.2}
\{(\nu_{i}\CS^{(\nu_{i})}_{\alpha}(\textbf{e});i\leq
N),\textbf{e}\}\overset{d}=\{(A^{(i)}_{g_{\textbf{e}}}; i\leq N),
L_{\textbf{e}}\},
\end{equation}
where $g_{t}\equiv\sup\{s<t:S_{s}=0\},$  which is equivalent to:

\begin{equation}\label{eq4.3}
\left\{ \begin{array}{rl}
 (i)& L_{\textbf{e}}=\textbf{e} \\
  (ii) & \{(A^{(i)}_{g_{\textbf{e}}};i\leq
  N)|L_{\textbf{e}}=l\}=\{(\nu_{i}\CS^{(\nu_{i})}_{\alpha}(\textbf{e});i\leq
N)\}
       \end{array} \right\}.
\end{equation}

In the next section, we give an explanation of~\mref{eq4.3}.
\section{An explanation of (4.3)}
We shall now show that (4.3) follows from the conjunction of the
next lemma and again the scaling property of a Bessel spider.
\subsection{Lemma}
\begin{lem}Let $(X_{t},t\ge 0)$ denote a nice Markov proces taking
values in a general space E; let $0$ belong to $E$, such that $0$
is regular for itself (with respect to $X$), and denote by
$(L_{t}, t\ge 0)$ a choice of its local time at 0; let
$(\tau_{\ell}, \ell\ge 0)$ denote the inverse local time, i.e.
$$
\tau_{\ell}=\inf\{t:L_{t}>\ell \}.
$$
It is a subordinator with associated Bernstein function,
$\psi(\theta)$:
$$
\E[{\mbox e}^{-\theta \tau_{\ell}}]={\mbox e}^{-\ell
\psi(\theta)}.
$$
Then, if $\mathbf{e_{\theta}}$ denotes an exponential $(\theta)$
random variable independent from $(X_{t},t\ge 0)$, and
$(A_{t},t\ge 0)$ is a continuous additive functional, one has
\begin{enumerate}
\item[(i)] $L_{\mathbf{e_{\theta}}}$ is exponentially distributed with
parameter $\psi(\theta);$
\item[(ii)] The conditional formula holds;
$$
\E[{\mbox
e}^{-A_{g_{\mathbf{e_{\theta}}}}}|L_{\mathbf{e_{\theta}}}=\ell]=\frac{E[{\mbox
e}^{-A_{\tau_{\ell}}-\theta \tau_{\ell}}]}{E[{\mbox e}^{-\theta
\tau_{\ell}}]}.
$$
\end{enumerate}
\end{lem}
This lemma is quite classical (it may be obtained from excursion
theory), let us admit it for a moment and derive (4.3) from it.
For this purpose, we take:
$A_{t}=\sum_{j=1}^{N}\lambda_{j}A^{(j)}_{t}, \lambda_{j}\ge 0,$
and we write (ii) in the lemma as follows:
$$
\E[\exp(-\sum_{j}\lambda_{j}A^{(j)}_{g_{\mathbf{e_{\theta}}}})|L_{\mathbf{e_{\theta}}}=l]
=\E[\exp(-\sum_{j}(\lambda_{j}+\theta)A^{(j)}_{\tau_{\ell}})]\exp(l\theta^{\alpha})
$$
(since: $\tau_{\ell}=\sum_{j}A^{(j)}_{\tau_{\ell}}$)
$$
=\exp\left(-l\sum_{j}p_{j}({(\lambda_{j}+\theta)}^{\alpha}-\theta^{\alpha})\right)
$$
(since the Ito excursion measure \textbf{n} associated with the
local time $(L_{t})$ satisfies
$\textbf{n}=\sum_{j}p_{j}\textbf{n}_{j},$ where $\textbf{n}_{j}$
is \textbf{n} restricted to excursions on $I_{j}$)
$$
=\E\left[\exp\left(-\sum_{j}\lambda_{j}\nu_{j}T^{(\nu_{j}\theta)}_{\alpha}(l)\right)\right]
$$
and it now suffices to take $\theta=1.$
\subsection{ A global view of the lemma}
In fact, the lemma is only a part of the more general statement.
\begin{prop} (We use the same notation as in the Lemma). The
following holds:
\begin{enumerate}
\item[(i)] The pre-$g_{\mathbf{e_{\theta}}}$ process $\{X_{u}, u\leq
g_{\mathbf{e_{\theta}}}\}$ and the post-$g_{\mathbf{e_{\theta}}}$
process $\{X_{g_{\mathbf{e_{\theta}}}+v}, v\leq
d_{\mathbf{e_{\theta}}}-g_{\mathbf{e_{\theta}}}\}$ where
$d_{t}=\inf\{s>t:X_{s}=0\},$ are independent;
\item[(ii)]The law of $(X_{u}, u\leq g_{\mathbf{e_{\theta}}})$ may
be described as follows:
 \item[(a)] $L_{\mathbf{e_{\theta}}}\equiv
 L_{g_{\mathbf{e_{\theta}}}}$ is $\exp(\psi(\theta))$ distributed;
 \item[(b)]$\{(X_{u}, u\leq
 g_{\mathbf{e_{\theta}}})|L_{\mathbf{e_{\theta}}}=l\}$ is
 distributed as: $(X_{u},u\leq \tau_{\ell})$ under the probability
 $$\exp(-\theta \tau_{\ell}+l\psi(\theta))\cdot P$$
 \item[(iii)]The law of the process $\{X_{g_{\mathbf{e_{\theta}}}+v}, v\leq
d_{\mathbf{e_{\theta}}}-g_{\mathbf{e_{\theta}}}\}$ is
$$
\frac{1}{\psi(\theta)}(1-\exp(-\theta
V(\varepsilon)))\textbf{n}(d\varepsilon)
$$
where V($\varepsilon$) denotes the lifetime of the generic
excursion $\varepsilon.$
\end{enumerate}
\end{prop}
For details of the proof and further references, the reader may
consult Salminen, Vallois and Yor~(2006).
 \vskip0.2in
\vskip0.2in \centerline{\Heading References} \vskip0.2in

\tenrm
\def\smc{\tensmc}
\def\sl{\tensl}
\def\bf{\tenbold}
\baselineskip0.15in

\Ref  \by Barlow, M., Pitman, J. and Yor, M. \yr 1989 \paper Une
extension multidimensionnelle de la loi de l'arc sinus. In
\textit{S\'eminaire de Probabilit\'es XXIII} (Azema, J., Meyer,
P.-A. and Yor, M., Eds.), 294--314, Lecture Notes in Mathematics
\textbf{1372}. Springer, Berlin\EndRef

\Ref \by Bertoin, J., Fujita, T., Roynette, B., and Yor, M. \yr
2006 \paper On a particular class of self-decomposable random
variables: the duration of a Bessel excursion straddling an
independent exponential time. To appear in Prob. Math. Stat
\EndRef

\Ref \by Bondesson, L. \yr 1992 \paper Generalized gamma
convolutions and related classes of distributions and densities.
Lecture Notes in Statistics, \textbf{76}. Springer-Verlag, New
York \EndRef

\Ref \by Chaumont, L. and Yor, M. \yr 2003 \book Exercises in
probability. A guided tour from measure theory to random
processes, via conditioning. Cambridge Series in Statistical and
Probabilistic Mathematics, 13 \publ Cambridge University Press,
Cambridge\EndRef

 \Ref \by    Cifarelli, D. M. and Melilli, E. \yr    2000
\paper Some new results for Dirichlet priors \jour \AnnStat \vol
28 \pages 1390-1413\EndRef

\Ref \by Devroye, L. \yr 1990 \paper A note on Linnik's
distribution \jour Statist. Probab. Lett. \vol 9 \pages
305-306\EndRef

\Ref \by Devroye, L. \yr 1996 \paper Random variate generation in
one line of code, in: 1996 Winter Simulation Conference
Proceedings, ed. J.M. Charnes, D.J. Morrice, D.T. Brunner and J.J.
Swain, pp. 265-272, ACM\EndRef

\Ref \by Donati-Martin, C., Roynette, B., Vallois, P., and Yor, M.
\yr 2005 \paper On constants related to the choice of the local
time et $0$, and the corresponding It\^o measure for Bessel
processes with dimension $d=2(1-\alpha),$ $(0<\alpha<1).$ To
appear in Studia Scientiarum Mathematicarum Hungarica \EndRef

\Ref \by Feller, W. \yr 1966 \paper Infinitely divisible
distributions and Bessel functions associated with random walks
\jour SIAM J. Appl. Math. \vol 14 \pages 864-875 \EndRef

\Ref\by Fujita, T. and Yor, M. \yr 2006 \paper An interpretation
of the results of the BFRY paper in terms of certain means of
Dirichlet Processes. Preprint\EndRef

\Ref \by Huillet, T. \yr 2000 \paper On Linnik's continuous-time
random walks \jour J. Phys. A \vol 33 \pages 2631-2652\EndRef

\Ref \by James, L. F.\yr 2006 \paper Gamma tilting calculus for
GGC and Dirichlet means with applications to Linnik processes and
occupation time laws for randomly skewed Bessel processes and
bridges. Preprint, arXiv:math.PR/0610218\EndRef

\Ref \by Kanter, M. \yr 1975 \paper Stable densities under total
variation inequalities \jour \AnnProb \vol 3 \pages 697-707
\EndRef

\Ref \by Lamperti, J. \yr 1958 \paper An occupation time theorem
for a class of stochastic processes \jour Trans. Amer. Math.
Soc.\vol 88 \pages 380-387\EndRef

\Ref \by Roynette, B. and Yor, M. \yr 2006 \paper Compl\'ement a
[BFRY]. In preparation \EndRef

 \Ref \by Salminen, P., Vallois, R, and Yor, M. \yr 2006
\paper On the excursion theory for linear diffusions. Submitted to
the Japanese Journal of Maths; a special volume for Prof. K.
It\^o. Preprint, arXiv:math.PR/0612687\EndRef

 \Ref \by Thorin, O.\yr 1977 \paper On the
infinite divisibility of the lognormal distribution \jour Scand.
Actuar. J. \vol 3 \pages 121-148 \EndRef

\Tabular{ll}
Lancelot F. James\\
The Hong Kong University of Science and Technology\\
Department of Information and Systems Management\\
Clear Water Bay, Kowloon\\
Hong Kong\\
\rm lancelot\at ust.hk\\
\EndTabular
\medskip

\Tabular{ll}
Marc Yor\\
Laboratoire de Probabilit\'es et Model\`es Al\'eatoires\\
Universit\'e Pierre et Marie Curie\\
Bo\^ite Courrier 188\\
4, Place Jussieu.\\
75252 Paris Cedex 05 France \EndTabular

\end{document}